\documentclass[11pt]{article}
\usepackage{amsfonts,latexsym}
\newcommand{\cleqn}{\setcounter{equation}{0}}
\newcommand{\clth}{\setcounter{theorem}{0}}
\newcommand {\sectionnew}[1]{\section{#1}\cleqn\clth}
\newcommand{\beq}{\begin{equation}}
\newcommand{\eeq}{\end{equation}}
\newcommand{\beqa}{\begin{eqnarray}}
\newcommand{\eeqa}{\end{eqnarray}}
\newcommand{\beaa}{\begin{eqnarray*}}
\newcommand{\eaa}{\end{eqnarray*}}
\newcommand{\nn}{\hfill\nonumber}
\newcommand{\ov}{\overline}
\newcommand \nc {\newcommand}
\nc \proof {\noindent {\em{Proof.\/ }}}
\nc \qed {$\Box$\hfill}
\newtheorem{theorem}{Theorem}[section]
\newtheorem{lemma}[theorem]{Lemma}
\newtheorem{proposition}[theorem]{Proposition}
\newtheorem{corollary}[theorem]{Corollary}
\newtheorem{definition}[theorem]{Definition}
\newtheorem{example}[theorem]{Example}
\newtheorem{remark}[theorem]{Remark}
\newtheorem{conjecture}[theorem]{Conjecture}
\newtheorem{question}[theorem]{Question}
\nc \bth[1] { \begin{theorem}\label{t#1} }
\nc \ble[1] { \begin{lemma}\label{l#1} }
\nc \bpr[1] { \begin{proposition}\label{p#1} }
\nc \bco[1] { \begin{corollary}\label{c#1} }
\nc \bde[1] { \begin{definition}\label{d#1}\rm }
\nc \bex[1] { \begin{example}\label{e#1}\rm }
\nc \bre[1] { \begin{remark}\label{r#1}\rm }
\nc \bcon[1] { \begin{conjecture}\label{con#1}\rm }
\nc \bque[1] { \begin{question}\label{que#1}\rm }
\nc {\eth} { \end{theorem} }
\nc {\ele} { \end{lemma} }
\nc {\epr} { \end{proposition} }
\nc {\eco} { \end{corollary} }
\nc {\ede} { \end{definition} }
\nc {\eex} { \end{example} }
\nc {\ere} { \end{remark} }
\nc {\econ} { \end{conjecture} }
\nc {\eque} { \end{question} }
\nc \eqref[1] {{\rm{(\ref{#1})}}}
\nc \thref[1]{Theorem \ref{t#1}}
\nc \leref[1]{Lemma \ref{l#1}}
\nc \prref[1]{Proposition \ref{p#1}}
\nc \coref[1]{Corollary \ref{c#1}}
\nc \deref[1]{Definition \ref{d#1}}
\nc \exref[1]{Example \ref{e#1}}
\nc \reref[1]{Remark \ref{r#1}}

\def \a {\alpha}
\def \b {\beta}
\def \g {\gamma}
\def \A {{\mathcal A}}
\def \B {{\mathcal B}}

\def \K {{\mathcal K}}
\def \M {{\mathcal M}}
\def \R {{\mathcal R}}

\def \Cset {{\mathbb C}}
\def \Zset {{\mathbb Z}}

\def \ord { {\mathrm{ord}} }

\def \spec { {\mathrm{Spec}} }

\def \ad { {\mathrm{ad}} }

\def \bb {b}

\nc \Gr {Gr}
\nc \GRN { \Gr^{(N)} }
\nc \GRA[1] { \Gr_A^{(#1)} }   
\nc \GRAN { \GRA{N} }
\nc \GrA[1] { \Gr_A(#1) }
\nc \GrAa { \GrA{\alpha} }
\nc \GRB[1] { \Gr_B^{(#1)} }   
\nc \GRBN { \GRB{N} }
\nc \GrB[1] { \Gr_B(#1) }
\nc \GrBb { \GrB{\beta} }
\nc \GRMB[1] { \Gr_{MB}^{(#1)} }   
\nc \GRMBN { \GRMB{N} }
\nc \GrMB[1] { \Gr_{MB}(#1) }
\nc \GrMBb { \GrMB{\beta} }

\begin{document}
\title{{\LARGE\bf{ General methods for constructing bispectral operators }}}
\author{
B.~Bakalov
\thanks{E-mail: bbakalov@fmi.uni-sofia.bg}
\quad
E.~Horozov
\thanks{E-mail: horozov@fmi.uni-sofia.bg}
\quad
M.~Yakimov
\thanks{E-mail: myakimov@fmi.uni-sofia.bg}
\\ \hfill\\ \normalsize \textit{
Department of Mathematics and Informatics, }\\
\normalsize \textit{ Sofia University, 5 J. Bourchier Blvd.,
Sofia 1126, Bulgaria }     }
\date{}
\maketitle
\begin{abstract}
We present methods for obtaining new solutions to the bispectral problem.
We achieve this by giving its abstract algebraic version
suitable for generalizations. All methods are illustrated by new
classes of bispectral operators.
\end{abstract}
\vspace{-11cm}
\begin{flushright}
{\tt{ q-alg/9605011 }}
\end{flushright}
\vspace{10cm}
\setcounter{section}{-1}
\sectionnew{Introduction}
The bispectral problem of J.~J.~Duistermaat and F.~A.~Gr\"unbaum \cite{DG}
consists of finding all {\em bispectral\/} ordinary differential  operators,
i.e.\ operators $L(x,\partial_x)$ having a family of eigenfunctions
$\psi(x,z)$, which are also eigenfunctions for  another differential  operator
$\Lambda(z,\partial_z)$ in the spectral parameter:
\beqa
&&L(x,\partial_x) \psi(x,z) = f(z) \psi(x,z),
\label{0.1}\hfill\\
&&\Lambda(z,\partial_z) \psi(x,z) = \theta(x) \psi(x,z).
\label{0.2}\hfill
\eeqa
In \cite{BHY2} we constructed large families of solutions to this problem,
generalizing all previously known results (cf.\ \cite{DG, Z, MZ, W, G}, etc.).
We obtained them as special {\em Darboux transformations\/} (called
``polynomial'')
 from the most obvious solutions -- the (generalized) {\em Bessel\/} and
{\em Airy\/} ones.
We recall that the Bessel (respectively -- Airy) operators have the form $L =
x^{-N} P(x\partial_x)$ (respectively $L = P(\partial_x) - x$) for some
polynomial $P$ of degree $N$.

For fixed $\psi$, all operators $L(x,\partial_x)$ for which $\psi$
is an eigenfunction form a commutative algebra $\A_\psi$, called the {\em
spectral algebra}. Then $\spec \A_\psi$ is an algebraic curve \cite{BC} -- the
{\em spectral curve}. The dimension of the space of eigenfunctions $\psi$
(= g.c.d.\ $\ord L$, $L\in\A_\psi$) is called
a {\em rank\/} of $\A_\psi$.
Following G.~Wilson \cite{W} we call the spectral algebra
$\A_\psi$ {\em bispectral\/} iff there exists an operator
$\Lambda(z,\partial_z)$ satisfying \eqref{0.2}.

The purpose of this paper is to present general methods for constructing
bispectral operators (respectively -- algebras). Some of our statements are
only abstract versions of concrete results
\cite{DG, W, BHY2, KR}. But adopting the algebraic point of view enables us to
clarify the ideas and to
give many classes of new examples of bispectral operators. Moreover, being
purely algebraic, some of the methods work also for difference or
$q$-difference analogs. The latter case is illustrated here by many examples
but more complete theory in the spirit of \cite{BHY2} will be considered
elsewhere.

In the process of our work \cite{BHY2} we substantially used the bispectral
involution $b$ introduced by G.~Wilson \cite{W}. It is, roughly speaking, the
map exchanging the roles of the arguments in the wave function $\psi(x,z)$:
\beq
b \psi(x,z) = \psi(z,x).
\eeq
 The main idea which we exploit here is to consider the bispectral involution
as an {\em anti-isomorphism\/} of algebras of ordinary differential operators
(cf.\ \cite{KR}), or more generally -- arbitrary rings.

 Our first construction of new bispectral algebras from old ones is
based on exhibiting anti-isomorphisms of algebras as compositions of already
known ones and inner {\em automorphisms}. In particular, we find new bispectral
algebras of arbitrary non-prime rank with both spectral curves (in $x$ and in
$z$) isomorphic to $\Cset$.

 The second construction is an abstract version -- which again uses only the
bispectral anti-isomorphism -- of the bispectrality theorem from \cite{BHY2}.
Analizing the latter has led us to the notion of {\em bispectral\/} Darboux
transformations which is applicable to other (than differential) operators.

 At the end we consider sequences of Darboux  transformations depending on
integers parametrizing the elements of the sequence. Then allowing the
parameters to be any complex numbers we again get new solutions to the
bispectral problem.


\sectionnew{Bispectral triples}
Here we present a general formulation of the bispectral problem which works not
only for differential but also for difference or $q$-difference modifications
of the original problem of \cite{DG}.

Let $\R$ and $\R'$ be associative algebras over a field $k$ and let $\M$ be a
left module over both of them. Let $\K$ and $\K'$ be fields such that
$k\subset\K\subset\R$ and $k\subset\K'\subset\R'$.
\bde{1.1}
We call an element $L\in\R$ {\em bispectral\/} iff there exist $\psi\in\M$,
$\Lambda\in\R'$, $f\in\K'$, $\theta\in\K$ such that
\beqa
L\psi=f\psi,
\label{1.1}\hfill\\
\Lambda\psi=\theta\psi.
\label{1.2}\hfill
\eeqa
\qed
\ede
In order to have non-trivial problem we assume that
\beq
L\psi=0 \;\; {\textrm{ implies }} \; L=0 \; {\textrm{ for }} L\in\R \;
({\textrm{or }} \R').
\label{1.3}
\eeq

Let us fix $\psi\in\M$ satisfying \eqref{1.1} and \eqref{1.2}. We are
interested in the equation (cf.\ \cite{KR})
\beq
P\psi=Q\psi
\label{1.4}
\eeq
for $P\in\R$, $Q\in\R'$. We put
\beqa
&&\B_\psi = \{ P\in\R \; | \; \exists Q\in\R' {\textrm{ for which \eqref{1.4}
is satisfied}} \},
\label{1.5}\hfill\\
&&\B'_\psi = \{ Q\in\R' \; | \; \exists P\in\R {\textrm{ for which \eqref{1.4}
is satisfied}} \}.
\label{1.5'}\hfill
\eeqa
We shall assume that the actions of $\R$ and $\R'$ on $\M$ commute and that
\eqref{1.3} holds. Then $\B_\psi$ and $\B'_\psi$ are associative algebras over
$k$ without zero divisors. Obviously \eqref{1.4} defines an anti-isomorphism
\beq
b\colon \B_\psi \to \B'_\psi, \quad b(P)=Q.
\label{1.6}
\eeq
Introduce also the subalgebras
\beqa
\K_\psi = \B_\psi \cap \K, \quad \K'_\psi = \B'_\psi \cap \K',
\label{1.7}\hfill\\
\A_\psi = b^{-1}(\K'_\psi), \quad \A'_\psi = b(\K_\psi).
\label{1.8}\hfill
\eeqa
Then $\A_\psi$ is a commutative algebra isomorphic to $\K'_\psi$ (the {\em
spectral algebra\/}). With a tilde we shall denote the fields of quotients of
these commutative rings, e.g.\ $\K_\psi \subset \widetilde\K_\psi \subset \K$.
Then $b$ extends to isomorphisms $\widetilde\K_\psi \to \widetilde\A'_\psi$ and
$\widetilde\A_\psi \to \widetilde\K'_\psi$.
\bde{1.2}
We call the triple $(\B_{\psi}, \B'_{\psi}, b)$ a {\em{bispectral
triple}} iff both $\K_{\psi}$ and $\K'_{\psi}$ contain non-zero
elements.
\qed
\ede
Obviously if $(\B_{\psi}, \B'_{\psi}, b)$ is a bispectral  triple then any
element $L \in \A_{\psi}$ is bispectral.
\bex{1.3}
(i) Let $\R \cong \R'$ be two copies of the rings of differential
operators of one variable with rational coefficients and let $\M$ be a module
of differentiable functions $\psi(x, z)$ (or formal power series) on which
$\R$ acts on the variable $x$ and $\R'$ acts on the variable $z$. Put
$\psi(x,z)= \delta (x-z)$ (Dirac delta function). Then $\B_\psi$ and $\B'_\psi$
are two copies of the Weyl algebra of the differential operators with
polynomial coefficients and $b$ is the formal conjugation $*$, i.e.\ $b(x)=z$,
$b(\partial_x)= -\partial_z$.

(ii) \cite{BHY2, KR}  Let $L_{\a}$ be the (generalized higher) Airy operator
$$
L_{\a}= \partial^N + \sum_{i=2}^{N-1} \a_i \partial^{N-i} - x,
\quad \a \in \Cset^{N-2}
$$
and let $\Phi(x)$ be a solution of the equation $L_\a\Phi=0$. Put $\psi(x,
z)=\Phi(x+z)$. Then $\B_\psi \cong \B'_\psi$
are again two copies of the Weyl algebra but now $b$ is
defined via
$$ b(x)=L_{\a}, \quad b(\partial)= \partial.$$
Then $\K_{\psi} = \Cset [x], \quad \A_{\psi}= \Cset [L_{\a}].$

(iii) \cite{BHY2} Let $L_{\b}$ be the (generalized) Bessel operator
$$L_{\b}= x^{-N} (D-\b_1) \cdots  (D-\b_N), \quad \b \in \Cset^N,$$
where $D=x\partial$. Let $\Phi(x)$ be a solution of the equation $L_{\b}
\Phi(x) = \Phi(x)$ and $\psi(x,z)=\Phi(xz).$ Then for ``generic'' $\b$ (see
\cite{BHY2}) $\B_\psi \cong \B'_\psi$ is generated by $x^N,D,L_\b$; now $b(x^N)
= L_\b$, $b(D)=D$, $b(L_\b) = x^N$ and $\K_\psi = \Cset[x^N]$, $\A_\psi =
\Cset[L_\b]$.
\qed
\eex
\bex{1.4}
(i) Let $\R$ and $\R'$ be two copies of the algebras of $q$-difference
operators with rational coefficients, i.e.\ operators of the type
$$D_q^N + a_1(x) D_q^{N-1} + \cdots + a_N(x),$$
where $D_q$ is the dilatation operator, $D_q f(x) = f(qx)$, $q \in \Cset$, and
$a_i(x)$ are rational functions. Let $\psi(xz) =\exp_q(xz)$ be the
$q$-exponent (see e.g.\ \cite{GR})
$$\exp_q(x) = \sum_{n=0}^\infty \frac{x^n}{(q;q)_n},$$
where $(q;q)_n=(1-q)(1-q^2) \cdots (1-q^n)$, $(q;q)_0 = 1$. Let $\B \cong \B'$
be two copies of the $q$-deformed Weyl algebra of operators
$$\sum a_j(x) \partial_q^j,$$
where $\partial_q=(D_q-1)/{x},$ $a_j(x) \in \Cset [x].$ Put $b(x)=-
\partial_q,$ $b(\partial_q)=-x$. Then $(\B, \B', b)$ is a bispectral
triple.

(ii) Let $L_{q,\b},$ $\b \in \Cset^N$, be the $q$-deformed Bessel operator:
$$L_{\b,q}= x^{-N} (D_q-q^{\b_1}) \cdots  (D_q-q^{\b_N}).$$
One can find solution $\psi(x)$ to the equation
$$L_{\b,q}\psi= \psi$$
of the form
$$\psi(x) = x^{\b_k} \sum_{j=0}^{\infty} a_j x^j, \quad a_0=1.$$
Define the algebra $\B$ to be generated by $D_q$, $x^N$ and $L_{\b,q}$ and
$\B' \cong \B$ (in $\B'$ the variable will be denoted by $z$). Then the
anti-izomorphism will be
    \beqa
           &&L_{\b,q}  \mapsto z^N       \hfill\\
 \bb:      &&x^N       \mapsto L_{\b,q}  \hfill\\
           &&D_q       \mapsto D_q      \hfill
    \eeqa
with a joint eigenfunction $\psi(x,z) = \psi(xz)$.
\qed
\eex
\bex{1.5}
Let $\psi(x,n)= 2^{-n} H_n(x),$ $x \in \Cset$, $n \in \Zset$, where
$H_n$ are the Hermite polynomials (see e.g.\ \cite{BE}) for $n \geq 0$ and
$H_n=0$ for
$n<0.$ Put $\B$ to be spanned by $x, \; \partial_x$ and $\B'$ to be spanned by
$n, \; T, \; T^{-1},$ where $T$ is the shift operator on functions defined on
$\Zset: \: (T f)(n) = f(n+1).$ Then an anti-isomorphism $b$ can be given by
    \beqa
  \bb:      &&x         \mapsto \left(T + \frac{n}{2} T^{-1}\right) \hfill\\
            &&\partial  \mapsto n T^{-1}.   \hfill
    \eeqa
It is easy to check that the Hermite operator $\partial^2 - 2x\partial$ is
bispectral. One can repeat this construction for any set of orthogonal
polynomials.
\qed
\eex

This example indicates that a general theory should work in various situations.
In the next section we will see how to build new bispectral triples from a
given one.
\sectionnew{Automorphisms of algebras of differential operators}

In this section we will assume that the algebras $\B$ and $\B'$ consist of
differential operators. The main tool for constructing new bispectral operators
will be the following simple observation.
\bpr{2.1}
Let $(\B, \B', \bb)$ be a bispectral triple and let $ L \in \B$.
If $\ad L$ is a locally nilpotent operator $\ad L : \B \rightarrow \B$ then
$\sigma = e^{\ad L}$ is an automorphism of $\B$ and $(\B, \B',
\bb \circ \sigma)$ is a bispectral triple.
\epr
(Here {\em locally nilpotent\/} means that for each $M\in\B$ there exists $n$
such that $(\ad L)^n M = 0$.)
The {\em proof\/} of \prref{2.1} being obvious is omitted.
\qed

\smallskip\noindent

An important example of an operator $L$ acting locally nilpotently on $\B$ is
any bispectral operator \cite{DG}. This can be proved by using the bispectral
anti-isomorphism $b$ and the fact that the function $f = \bb(L)$ acts
locally nilpotently on $\B'$.
\bre{2.1,5}
(i) Obviously, exchanging the roles of $\B$ and $\B'$ one can use also
the anti-isomorphism $\bb_2 = e^{\ad L'} \circ b$, where $L'$ acts
locally nilpotently on $\B'$. We will use \prref{2.1} in both versions.

(ii) In the definition of a bispectral triple a decisive role is played by the
function $\psi(x,z) \in \M$ via which the anti-isomorphism is defined.
\prref{2.1} defines a new anti-isomorphism $\bb_1$ in terms of the old
function $\psi$ (see \eqref{1.4}). But to have really a bispectral triple in
the sense of  \deref{1.1} we face the difficult task to define a new function
$\tilde \psi$ to satisfy $ P \tilde \psi = \bb_1(P) \tilde \psi$. The most
natural choice $\tilde \psi = e^L \psi$ often makes no sense. In the examples
below after defining the bispectral triple using $\psi$ the new function
(or formal power series) $\tilde \psi$ can be determined using other arguments.
\qed
\ere
Although trivial the above Proposition is a source of new bispectral
operators.
\bex{2.2}
(i) Let $\B$ and $\B'$ be two copies of the Weyl algebra and let
$\bb(x) = \partial_z,$ $\bb (\partial_x) = z,$ $\psi=e^{xz}$.
For any polynomial $p(\xi)$ the operators $\ad p(x)$ and
$\ad p(\partial_x)$ act locally nilpotently on $\B$. Then
if we put $\sigma_1 = e^{\ad p(x)}$ and $\bb_1=\bb \circ \sigma_1$
we have
    \beqa
&& \bb_1 (x) =b(x) = \partial_z,  \hfill\\
&& \bb_1 (\partial_x)  =
\bb(-p'(x)+\partial_x)=z-p'(\partial_z). \hfill
    \eeqa
In this way we get a new bispectral triple $(\B, \B', \bb_1)$ with a new
function $\tilde \psi (x,z) = e^{p'(x)} e^{xz}$.

(ii) One can apply the automorphism $\sigma_2 = e^{\ad q(\partial_x)}$ to the
triple
$(\B, \B', \bb_1)$. Then we get $\bb_2 = \bb_1 \circ \sigma_2$, which gives
    \beqa
&& \bb_2 (x)
=\bb_1(x+q'(\partial_x))=\partial_z+q'(z-p'(\partial_z), \hfill\\
&& \bb_2 (\partial_x)  = \bb_1 (\partial_x) = z-p'(\partial_z).\hfill
    \eeqa
For the new function $\psi_2(x,z)$ we must have
    \beqa
&&\partial_x \psi_2 = (z-p'(\partial_z)) \psi_2,
\label{2.1}\hfill\\
&&  x \psi_2 = (\partial_z-q'(z-p'(\partial_z))) \psi_2.
\label{2.2}\hfill
    \eeqa
This implies
    \beq
(\partial_x+p'(x-q'(\partial_x))) \psi_2 = z\psi_2,
\label{2.3}
    \eeq
which gives a new non-trivial example of bispectral operators, provided
$\psi_2$ exists. One can construct a solution $\psi_2$ using Fourier
transformation with respect to both $z$ and $x$ on \eqref{2.1} and on
    \beq
x \psi_2 = (\partial_z+q'(\partial_x)) \psi_2.
\label{2.4}
    \eeq
Because $\K_{\psi_2} = \Cset[x]$, $\K'_{\psi_2} = \Cset[z]$ we see that both
spectral curves are isomorphic to $\Cset$.
\qed
\eex
\bre{}
The automorphisms considered in \exref{2.2} generate the entire group of
automorphisms of the Weyl algebra \cite{D}.
This result leads to a classification of all solutions to the bispectral
problem such that $\B_\psi \cong \B'_\psi$ is the Weyl algebra.

In general, we propose the following approach to the classification problem.
Given solution $\psi$ to the bispectral problem is first reduced to a one with
a ``better'' $\B_\psi$ (perhaps by a ``bispectral Darboux
transformation'' -- see the next section). In all known to us cases this
$\B_\psi$ is either the Weyl algebra or $\langle x^N, D, L_\b
\rangle$ (see \exref{1.3}). Then the second step is for a given ring $\B$ to
find all bispectral triples $(\B,\B',b)$.
\qed
\ere
One can easily check that the operators in the r.h.s.\ of \eqref{2.1} and
\eqref{2.2} satisfy the so-called string equation $[P,Q] = 1$, which plays an
important role in $2$d-quantum gravity \cite{BK, DS, GM} (see also \cite{Wit,
K}). The corresponding
to the r.h.s.\ of \eqref{2.1} tau-functions can be shown to satisfy the
$W$-constraints (see \cite{AvM, FKN}). This will be done in another place.

An anti-isomorphism which can successfully be used for constructing new
bispectral operators from old ones is the formal conjugation: $\bb(x)=x$,
$\bb(\partial_x)=-\partial_x$. Instead of formulating a general result we shall
give an example.
\bex{2.3}
Let $\b \in \Cset^N, \; \g \in \Cset^N, \; \b_i+\g_i=N-1$.
Define $\B=\langle x^N, L_\b, D_x \rangle,$
$\B'=\langle z^N, L_\gamma, D_z \rangle,$ where the operators $L_\b, \; L_\g$
 are Bessel operators (see \cite{BHY2} or \exref{1.3} (iii)). Let $b$ be the
formal conjugation: $\bb(x^N) =z^N,$
$\bb(D_x) =-D_z-1,$ $\bb(L_\b) =(-1)^N L_\g.$ As $L_\g$ acts locally
nilpotently on $\B'$ (since it is bispectral) one can define $\bb_1=
e^{-\ad N^{-1} L_\g} \circ b$. Then we have
    \beqa
&& \bb_1(D_x)=-D_z-1 + L_\g, \nn \hfill\\
&& \bb_1(L_\b)= (-1)^N L_\g, \nn \hfill\\
&& \bb_1(x^N)=M_\g= e^{-\ad N^{-1} L_\g} (z^N).
      \nn \hfill
    \eeqa
$M_\g$ can be computed easily applying $b_1$ to both sides of the identity
$$
L_\b x^N = (D_x+N-\b_1) \cdots (D_x+N-\b_N).
$$
In the simplest nontrivial case $N=3$ we have
$$
M_\g = -L_\g^2 + \left(3 D_z -9 - \sum \g_i\right) L_\g
   - 3 D_z^2 + \left(9 + 2\sum \g_i\right) D_z - \left(9 + 3\sum\g_i + \sum
\g_i\g_j\right) + z^3.
$$
The corresponding function $\psi(x,z)$ satisfies
    \beqa
&& (D_x + D_z +1) \psi = L_\g \psi = (-1)^N L_\b \psi,
\nn \hfill\\
&& x^N \psi = M_\g \psi, \nn \hfill\\
&& z^N \psi = M_\b \psi. \nn \hfill
    \eeqa
For ``generic'' $\b$ (see \cite{BHY2}) both spectral curves are $\Cset$.

This construction can be generalized considering $b_1 = e^{\ad p(L_\g)} \circ
b$ for arbitrary polynomial $p$.
 \qed
\eex
\bex{2.6}
The last example deals with the non-generic $\b=(-1,1,3)$. In this case the
algebra $\B_\psi = \langle x^2, x^3, D_x, L_\b \rangle$ is larger than for
generic $\b$ (cf.\ \exref{1.3}~(ii)). We have the anti-isomorphism
\beqa
       && x^2 \mapsto z^2, \quad x^3 \mapsto z^3, \nn\\
b:     && D_x \mapsto -D_z-1, \nn\\
       && L_\b \mapsto L_\b. \nn
\eeqa
Let $b_1 = e^{(-1/3) \ad L_\b} \circ b$. We obtain the bispectral operator
$$
L = e^{(-1/3) \ad L_\b}(x^2) =
x^{-4}(D+1)(D-1)(D-2)(D-4) - 2x^{-1}(D+1)(D-1) + x^2,
$$
which can be written as
$$
L = (\partial^2 - x^{-1} \partial - x) (\partial^2 + x^{-1} \partial - x).
$$
It coincides with a Darboux transformation of the Airy operator $\partial^2
- x$ -- see \cite{BHY2}, Example~5.8 for $d_0=1, a=\infty, \lambda=\mu=0$ (cf.\
also \cite{G, KR}).

The spectral algebra has rank 2, the spectral curve is $\Cset$ with a cusp at
0.
\qed
\eex
We hope the above examples help in unifying the Bessel and Airy cases of
\cite{BHY2} which although similar still had important differences.

\sectionnew{Bispectral Darboux transformations}
Up to now the most widely used method for constructing solutions to the
bispectral
problem is the method of Darboux trasformations introduced in the pioneering
work of Duistermaat and Gr\"unbaum \cite{DG} (see also \cite{Z, MZ, BHY1,
BHY2, KR}).
Here we give a general definition of  bispectral Darboux transformations in the
context of bispectal triples introduced in Sect.~1.
\bde{3.1}
Let $(\B_{\psi}, \B'_{\psi} , \bb)$ be a bispectral triple.
Let $L \in \A_\psi$ be a
bispectal operator which can be factorized as
\beq
L=Q \theta^{-1} P \;{\textrm{ with }} \quad P, \, Q \in \B_\psi, \;
\theta \in \K_\psi.
\label{3.1}
\eeq
We call the operator
\beq
\ov L=P Q \theta^{-1}
\label{3.2}
\eeq
a {\em{bispectral Darboux transformation}} of $L$ and the function
\beq
\ov \psi = P \psi
\label{3.3}
\eeq
a bispectral Darboux transformation of $\psi$.
\qed
\ede
\bth{3.2}
The operator $\ov L$ is bispectral. More precisely, if
$f=\bb(L)\in\K'_\psi$ then
\beq
\ov L \, \ov \psi = f \, \ov \psi
\label{3.4}
\eeq
and
\beq
\ov \Lambda \, \ov \psi = \theta \, \ov \psi
\label{3.5},
\eeq
where
\beq
\ov \Lambda = \bb(P) \bb(Q) f^{-1}.
\label{3.6}
\eeq
\eth
\proof Let $\Lambda =\bb(\theta).$ As $\B_\psi$ has no zero divisors
\eqref{3.1} implies
$$ \theta = P L^{-1} Q, \quad L^{-1} \in \widetilde \A_\psi .$$
Applying the anti-isomorphism $\bb$ we obtain
$$ \Lambda =\bb(\theta)= \bb(Q) f^{-1} \bb(P).$$
This shows that $\ov \Lambda$ is a bispectral Darboux transformation of
$\Lambda$. From \eqref{3.2} and \eqref{3.6} it follows that $\ov \psi$
satisfies \eqref{3.4} and \eqref{3.5}.
\qed

\smallskip\noindent

\thref{3.2} is an abstract version of our bispectrality theorem from
\cite{BHY2}.
However, it formalizes only partially the arguments of this theorem. The main
difference is that here we do not touch the subtle (in general) question
when a factorization of the form \eqref{3.1} is possible. Nevertheless,
\thref{3.2} can be a part of a general theory of bispectral operators. Below we
illustrate it by examples which have not appeared elsewhere.
\bex{3.3}
(i) Let $(\B, \B', \bb)$ be the bispectral triple from \exref{1.4}~(i), i.e.\
$\B \cong \B'$ be two copies of the $q$-deformed Weyl algebra. Then
$L= \partial_q^2$ is a bispectal operator with $\psi(x,z) = \exp_q(xz)$, i.e.\
$$
\partial_q^2 \exp_q(xz) = z^2 \exp_q(xz).
$$
One can easily check that
$$
L = \bigl( (a + q^2 x) \partial_q + q^2 - 1 \bigr)
\bigl( (a + x)(a + q x) \bigr)^{-1}
\bigl( (a + x) \partial_q - q + 1 \bigr),
$$
which gives the factorization of \eqref{3.1}. Then if we put
$$
\psi_1(x,z)= ((a+x) \partial_q - q + 1) \exp_q (xz),
$$
\thref{3.2} gives that
    \beqa
&&\bigl( (a+x)\partial_q -q +1 \bigr)
\bigl( (a + q^2 x) \partial_q + q^2 - 1 \bigr)
\bigl( (a + x)(a + q x) \bigr)^{-1}
\psi_1(xz) = z^2 \psi_1 (xz),
  \nn\hfill\\
&&\bigl( -z (a-\partial_q) -q +1 \bigr)
\bigl( -z (a-q^2 \partial_q) + q^2 -1 \bigr)
z^{-2} \psi_1 (xz) = (a+x)(a+qx) \psi_1(xz).
  \nn\hfill
     \eeqa

(ii) Take the same bispectral triple but put
$$
L = \partial_q^3 = x^{-3} (D-1)(q^{-1}D-1)(q^{-2}D-1).
$$
For $\{ a,b,c \} = \{ 0,1,2 \}$
$L$ can be factorized as
$(q^{-a+3} D-1) (q^{-b+3} D-1) x^{-3} (q^{-c} D-1)$.
The procedure of \thref{3.2} gives the operator
$$
\ov L = x^{-3} (q^{-a} D-1) (q^{-b} D-1) (q^{-c-3} D-1).
$$
It is clear that repeatedly applying this procedure  one can get
$$
\widetilde L = x^{-3} ( q^{-3k}D-1)(q^{-3m-1}D-1) (q^{-3n-2}D-1)
$$
with any $k, m, n \geq 0.$ Notice that $\widetilde L$ is a $q$-Bessel operator
defined in \exref{1.4}~(ii).

(iii) Let $(\B, \B', b)$ be the bispectral triple defined in \exref{1.4}~(ii)
 and put
$$
L = L_{\b,q} = x^{-2} (D - q^{\b_1})(D - q^{\b_2}).
$$
Let $\b_2 - \b_1 = 2 \a$, $\a\in\Zset_{\ge0}$.

Then one can check that $L$ factorizes as
$
L = Q \theta^{-1} P
$
with
\beqa
Q &=& q^2 D (1+a q^{2\a} x^{2\a}) - q^{\b_2} (1 + a x^{2\a}),
\nn\\
P &=& (1+a x^{2\a}) D - q^{\b_1} (1 + a q^{2\a} x^{2\a}),
\nn\\
\theta &=& x^2 (1+a x^{2\a}) (1+a q^{2\a} x^{2\a}),
\nn
\eeqa
for $a\in\Cset$.
Applying the recipe of \thref{3.2} one gets easily the bispectral operators
$\ov L, \; \ov \Lambda :$
$$
\ov L = P Q \theta^{-1}, \qquad \ov\Lambda = b(P) b(Q) f^{-1},
$$
satisfying (\ref{3.4}, \ref{3.5}), where $f=z^2$ and
\beqa
b(P) &=& D (1+a L^\a) - q^{\b_1} (1 + a q^{2\a} L^\a),
\nn\\
b(Q) &=& (1+a q^{2\a} L^\a) q^2 D - q^{\b_2} (1 + a L^\a).
\eeqa
\qed
\eex

We end this section with the following observation. Suppose starting with a
bispectral operator one produces a sequence of bispectral Darboux
transformations labeled by some integer parameters $n_1, \ldots , n_k$.
Suppose that all the operators $L_{n_1, \ldots , n_k}$ are of the same order
and that their coefficients depend rationally on $n_1, \ldots , n_k$. Then
considering $n_1, \ldots , n_k$ as arbitrary complex parameters we get a family
of bispectral operators. The standard example is from \cite{DG} where from
$\partial^2$ one produces $\partial^2 - n(n-1)x^{-2}$ and putting
$n(n-1)=c$ one gets arbitrary classical Bessel operator. Another example is
provided by \exref{3.3}~(ii). The following is a non-trivial example
illustrating this method.
\bex{3.4}
If $\psi(x, z)$ is the ordinary Airy solution of the bispectral
problem, i.e.\
    \beqa
&&(\partial_x^2 -x) \psi = z \psi, \nn\hfill\\
&&\partial_x \psi = \partial_z \psi, \nn\hfill\\
&& x \psi = (\partial_z^2-z) \psi, \nn\hfill
    \eeqa
then
$D_z \psi = \partial_x ( \partial_x^2-x) \psi =
(x^{-3} D_x (D_x -1)(D_x - 2) - D_x -1) \psi.$
{}From this it is easy to see that $\psi$ is the joint
eigenfunction of \exref{2.3} for $\b=(0,1,2)$. For arbitrary $\b \in
\Cset^3$
denote that function by $\psi_\b$. Then we have the following Darboux
transformation:
$$ \psi_{\b_1+1, \b_2+1, \b_3-2}=
   \frac{1}{z} \biggl( L_{\b_1,\b_2}(x, \partial_x) -x \biggr)
   \psi_{\b_1, \b_2, \b_3}. $$
For example the identity
$$ D_z \psi_{\b_1+1, \b_2+1, \b_3-2}=
   \biggl( L_{\b_1+1,\b_2+1,\b_3-2}(x, \partial_x) - D_x -1 \biggr)
   \psi_{\b_1+1, \b_2+1, \b_3-2} $$
is a straightforward consequence of the equation
$$
   (L_{\b_1,\b_2}-x) (L_{\b_1,\b_2,\b_3}-D-2)=
   (L_{\b_1+1,\b_2+1,\b_3-2}-D-1)(L_{\b_1,\b_2}-x).
$$
In this way we obtain the bispectrality of \exref{2.3} for
$\b=(i, 1+j, 2-i-j)$, where $i, j$ are arbitrary integers.
Because of the polynomial dependence on $\b$ of the above differential
operators we can continue the joint eigenfunction for arbitrary $\b \in
\Cset^3$ such that $\b_1+\b_2+\b_3 = 3$.
Note that in \exref{2.3} one can
normalize $\b_1+\b_2+\b_3$ by
multiplying $\psi_\b$ by $x^c z^{-c}$ for $c \in \Cset$ (cf.\ \cite{BHY2}).
\qed
\eex
The last Example together with Examles \ref{e2.3} and \ref{e2.6} shows that
there is a connection among all methods presented here.

{\flushleft{\bf{Acknowledgement}}}

\medskip\noindent
We thank Alex Kasman for stimulating correspondence and for sending us his
preprint \cite{KR}. This work was partially supported by Grant MM-523/95
of Bulgarian Ministry of Education, Science and Technologies.
\begin{small}
    
\end{small}

\end{document}